\documentclass{amsart}
\usepackage{amssymb}

\usepackage{amsfonts}
\usepackage{amssymb}
\newtheorem{theorem}{Theorem}[section]

\newtheorem{claim}[theorem]{Claim}
\newtheorem{conjecture}[theorem]{Conjecture}

\newtheorem{definition}[theorem]{Definition}
\newtheorem{question}[theorem]{Question}
\newtheorem{remark}[theorem]{Remark}

\newcommand{\BC}{\mathbb{C}}
\newcommand{\BP}{\mathbb{P}}
\newcommand{\BQ}{\mathbb{Q}}

\newcommand{\BZ}{\mathbb{Z}}

\newcommand{\Aut}{\operatorname{Aut}}

\newcommand{\id}{\operatorname{id}}
\newcommand{\Imm}{\operatorname{Im}}

\newcommand{\NS}{\operatorname{NS}}
\newcommand{\ord}{\operatorname{ord}}
\newcommand{\Pic}{\operatorname{Pic}}
\newcommand{\ttop}{\operatorname{top}}

\newcommand{\SC}{\mathcal{C}}

\newcommand{\Alb}{\operatorname{Alb}} 
\newcommand{\OH}{\operatorname{H}} 
\newcommand{\LS}[1]{|#1|} 

\newcommand{\ratmap}
{{\,\cdot\negmedspace\cdot\negmedspace\cdot\negmedspace\to\,}}

\newcommand{\alg}{\mathrm{alg}}

\begin{document}
\begin{large}
\title
[Dynamics of Automorphisms of Compact Complex Manifolds] 
{Dynamics of Automorphisms of Compact Complex Manifolds}
\author{De-Qi Zhang}
\address
{
\textsc{Department of Mathematics} \endgraf
\textsc{National University of Singapore, 2 Science Drive 2,
Singapore 117543, Singapore
\endgraf
}}
\email{matzdq@nus.edu.sg}

%

\begin{abstract}
We give an algebro-geometric approach towards the dynamics of
automorphisms/endomorphisms of projective varieties or compact K\"ahler manifolds,
try to determine the building blocks of automorphisms /endomorphisms, and
show the relation between the dynamics of automorphisms/endomorphisms and the geometry
of the underlying manifolds.
\end{abstract}
\subjclass{14J50, 14E20, 32H50} 
\keywords{Automorphism, endomorphism, dynamics, topological entropy,
Calabi-Yau manifold,
rational connected manifold.
}
\maketitle
\section{Introduction}

%
%
%
%
%
%

\par
\S 1. Introduction

\par
\S 2. Preliminaries: entropy, dynamical degree

\par
\S 3. Gizatullin-Harbourne-McMullen conjecture

\par
\S 4. Tits alternative for automorphism group

\par
\S 5. Dynamics of automorphisms

\par
\S 6. Building blocks of endomorphisms

\par
\S 7. Cohomologically hyperbolic endomorphism

\par
References

\vskip 12mm

\section{Introduction} 

I will report recent results obtained in 
\cite{KOZ}, \cite{NZ}, \cite{Zh1}, \cite{Zh2}, \cite{Zh3},
some being jointly with J. Keum, N. Nakayama and K. Oguiso.
\par
Readers may first skip \S 2 and go directly to subsequent sections.
For entropy and dynamical degree, there is a short survey \cite{Fr06}.
For (dynamics of polarized) endomorphisms, please see surveys \cite{Zs} and \cite{FN}.

\par
We work over the field $\BC$ of complex numbers.

\newpage

\section{Preliminaries: entropy, dynamical degree}\label{pre} 

\par
We use the convention of Hartshorne's book, \cite{KMM} and \cite{KM}.

\par
Let $X$ be a compact complex K\"ahler manifold with
$$H^*(X, {\BC}) = \bigoplus\nolimits_{i\ge0}
H^i(X, \BC)$$ 
the total cohomology group.
Take a $g \in \Aut(X)$. Denote by $\rho(g)$ the {\it spectral
radius} of $g^*| H^*(X, \BC)$. It is known that either $\rho(g) > 1$,
or $\rho(g) = 1$ and all eigenvalues of $g^* | H^*(X, \BC)$ are of
modulus $1$. When $\log \rho(g) > 0$ (resp. $\log \rho(g) = 0$) we
say that $g$ is of {\it positive entropy} (resp. {\it null
entropy}).
\par
We refer to Gromov \cite{Gr}, Yomdin \cite{Yo}, Friedland
\cite{Fr95}, and Dinh - Sibony \cite{DS05}, for the
definition of the $i$-th {\it dynamical degree} $d_i(g)$ for $1 \le
i \le n = \dim X$ (note that $d_n(g) = 1$ now and we set $d_0(g) = 1$)
and the actual definition of the {\it topological entropy} $h(g)$
which turns out to be $\log \rho(g)$ in the current setting.
\par
\begin{remark}
{\rm The above terminlogies can also
be defined for endomorphisms and even for dominant meromorphic self-maps.}
\end{remark}
\par
Let $Y$ be a projective variety and $g \in \Aut(Y)$. We say that
$g$ is of {\it positive entropy}, or {\it null entropy}, or {\it
parabolic}, or {\it periodic}, or {\it rigidly parabolic}, or of
{\it primitively positive entropy} (see the definitions below), if
so is $g \in \Aut(\tilde{Y})$, where $\tilde{Y} \rightarrow Y$ is
one (and hence all) $g$-equivariant resolutions guaranteed by
Hironaka. The definitions do not depend on the choice of
$\tilde{Y}$ because every two $g$-equivariant resolutions are
birationally dominated by a third one.
\par
We use $g | Y$ to signify that $g \in \Aut(Y)$.
\par
By a {\it pair} $(Y, g)$ we mean a projective
variety $Y$ and an automorphism $g \in \Aut(Y)$. Two pairs $(Y', g)$
and $(Y'', g)$ are {\it birationally equivariant}, if there is
a birational map $\sigma : Y' \cdots \to Y''$ such that
the action $\sigma (g|Y') \sigma^{-1} : Y'' \cdots \to Y''$
is biregular.
\par
$g \in \Aut(Y)$ is {\it periodic} if the order $\ord(g)$ is
finite. $g$ is {\it parabolic} if $\ord(g) = \infty$ and if $g$ is
of null entropy.
\par
$(Y', g)$ is {\it rigidly parabolic} if ($g | Y'$ is parabolic
and) for every pair $(Y, g)$ which is birationally equivariant
to $(Y', g)$ and for every $g$- equivariant surjective morphism $Y
\rightarrow Z$ with $\dim Z > 0$, we have $g | Z$ parabolic.
\par
Let $Y'$ be a projective variety and $g \in \Aut(Y')$ of
positive entropy (so $\dim Y' \ge 2$). A pair $(Y', g)$ is of {\it
primitively positive entropy} if it is not of imprimitive positive
entropy, while a pair $(Y', g)$ is of {\it imprimitively positive
entropy} if it is birationally equivariant to a pair $(Y, g)$
and if there is a $g$-equivariant surjective morphism $f: Y
\rightarrow Z$ such that either one of the two cases below occurs.
%
%
\begin{itemize}
\item[(a)]
$0 < \dim Z < \dim Y$, and $g | Z$ is still of positive entropy.
\item[(b)]
$0 < \dim Z < \dim Y$, and $g | Z$ is periodic.
\end{itemize}
\par
\begin{remark}
{\rm
We observe that in Case(b), for some $s
> 0$ we have $g^s | Z = \id$ and that $g^s$ acts faithfully on the
general fibre $Y_z$ of $Y \rightarrow Z$, such that $g^s | Y_z$ is
of positive entropy; see \cite[(2.1)(11)]{Zh2}.
In fact, we have $d_1(g^s|Y) = d_1(g^s|Y_z)$;
see \cite[Appendix, Theorem A.10]{NZ}.
}
\end{remark}

\par
\begin{remark}
{\rm
In view of the observation above,
if $\dim Y \le 2$ and if the pair $(Y, g)$ is of
positive entropy, then $\dim Y = 2$ and the pair $(Y, g)$ is always
of primitively positive entropy.
}
\end{remark}

\par
For the references to the following important results,
see Dinh-Sibony \cite{DS05} and Guedj \cite{Gu05}.

\begin{theorem} \label{theorem:ent}
Let $M$ be a compact K\"ahler manifold of dimension $n$ and let
$g$ be an automorphism of $M$. By $\rho(g^{*} \vert W)$ we denote
the spectral radius of the action of $g^*$ on a $g^*$-stable
subspace $W$ of the total cohomology group $H^{*}(M, \BC)$.
Then we have\emph:

\begin{itemize}
\item[$(1)$] $\rho(g^{*} \vert H^{*}(M, \BC)) \ge 1$, and $\rho(g^{*}
\vert H^{*}$ $(M, \BC)) = 1$ (resp. $> 1$) if and only if
$\rho(g^{*} \vert H^{2}(M, \BC)) = 1$ (resp. $> 1)$.
Further, $\rho(g^{*} \vert H^{*}(M, \BC))$ $= 1$
(resp. $> 1$) if and only if so is for $g^{-1}$.

\item[$(2)$] One has 
$$\rho(g^{*} \vert H^{2}(M, \BC)) = \rho(g^{*} \vert H^{1,1}(M)).$$
Further, if $M$ is projective, then the above value equals
$\rho(g^{*} \vert \NS(M))$.

\item[$(3)$] The spectral radius $\rho(g^*|H^{i,i}(X, {\BC}))$ equals
the $i$-th dynamical degree $d_i(g)$. Further, there are integers
$m \le m'$ such that
$$1 = d_0(g) < \dots < d_m(g) = \dots = d_{m'}(g)
> \dots > d_n(g) = 1.$$

\item[$(4)$] The following holds with $h(g)$ the topological entropy of $g$\emph:
$$\rho(g^* | H^*(X, {\BC})) = \max_{0 \le i \le n} \, d_i(g) = e^{h(g)}.$$
\end{itemize}
\end{theorem}

\section{Gizatullin-Harbourne-McMullen conjecture}

%
%
%
%

Consider the following question of Gizatullin-Harbourne-McMullen
(see \cite[page 409]{Ha} and \cite[\S 12]{Mc07}),
where $\Aut^*(X) := \Imm(\Aut(X) \to \Aut(\Pic(X)))$:

\begin{question}\label{QGHM}
{\rm
Let $X$ be a smooth projective rational surface.
If $\Aut^*(X)$ is infinite, is there then a birational morphism $\varphi$ of $X$ to
a surface $Y$ having an anti-pluricanonical curve and an infinite subgroup $G \subset \Aut^*(Y)$
such that $G$ lifts via $\varphi$ to $X$?
}
\end{question}

\begin{remark}
{\rm
It is very difficult to construct examples of automorphisms $g$ on rational surfaces $X$
with positive entropy. There are some sporadic examples by A. Coble and M. Gizatullin.
In \cite{Mc07}, McMullen succeeded in constructing an infinite series of such examples
$(X_n, g_n)$. In all the examples of his, $X_n$ has an anti-canonical curve.
This lead him to ask a question similar to the one above.
In \cite{Zh1}, we show that $\Aut(X_n) = (\text{\rm finite group}) \rtimes \BZ$,
so each of McMullen's surfaces essentially has only one automorphism of positive entropy.
}
\end{remark}

\par
A member in an anti-pluricanonical system
$|{-}nK_X|$ ($n\ge 1$) is called an {\it anti n-canonical} curve (or {\it divisor})
or an {\it anti-pluricanonical} curve;
a member in $|{-}K_X|$ is an {\it anti $1$-canonical} curve,
or an {\it anti-canonical} curve.
\par
The result below answers Question \ref{QGHM}
in the case of null entropy.
\par \vskip 0.5pc
\begin{theorem}[{cf. \cite{Zh1}}]\label{Thnull}
Let $X$ be a smooth projective rational surface and $G \le \Aut(X)$
an infinite subgroup of null entropy.
Then we have\emph:
\begin{itemize}
\item[$(1)$]
There is a $G$-equivariant smooth blowdown $X \rightarrow Y$
such that $K_Y^2 \ge 0$ and hence $Y$ has an anti-pluricanonical curve.
\item[$(2)$]
Suppose further that $\Imm (G \rightarrow \Aut(\Pic(X)))$
is also an infinite group. Then the $Y$ in $(1)$ can be so chosen that
$-K_Y$ is nef of self intersection zero
and $Y$ has an anti $1$-canonical curve.
\end{itemize}
\end{theorem}
\par \vskip 0.5pc
For groups which are not necessarily of null entropy,
we have the following result
which is especially applicable (with the same kind of $H$)
when $G/H \ge \BZ \rtimes \BZ$.
%
%
\par \vskip 0.5pc
\begin{theorem}[{cf. \cite{Zh1}}] \label{Thz2}
Let $X$ be a smooth projective surface and $G \le \Aut(X)$ a subgroup.
Assume that there is a sequence of groups
$$H \lhd A \lhd G$$
satisfying the following three conditions\emph:
\begin{itemize}
\item[$(1)$]
$\Imm(H \rightarrow \Aut(\NS(X)))$ is finite;
\item[$(2)$]
$A/H$ is infinite and abelian; and
\item[$(3)$]
$|G/A| = \infty$.
\end{itemize}
Then $G$ contains a subgroup $S$ of null entropy and infinite order.
\newline
In particular, when $X$ is rational, there is an $S$-equivariant
smooth blowdown $X \rightarrow Y$ such that $Y$ has an anti-pluricanonical curve.
\end{theorem}
\par \vskip 0.5pc

\newpage

\begin{remark}
$ $
{\rm
\begin{enumerate}
\item[$(1)$] 
Conditions like the ones in Theorems \ref{Thnull} and \ref{Thz2} are
necessary in order to have an affirmative answer to
Question \ref{QGHM}. See Bedford-Kim
\cite[Theorem 3.2]{BK} for a pair $(X, g)$
with $g$ of positive entropy and Iitaka D-dimension $\kappa(X, {-}K_X) = -\infty$.
\item[$(2)$] 
The blowdown process $X \rightarrow Y$ to the minimal pair $(Y, S)$
in Theorems \ref{Thnull} and \ref{Thz2} is necessary,
as observed by Harbourne \cite{Ha}.
\end{enumerate}
}
\end{remark}

\section{Tits alternative for automorphism group}
%
%
%
%

We often study a group $G$ of automorphisms of a projective
variety $X$ through its action on cohomological spaces such as on
the N\'eron-Severi group (over $\BQ$) $\NS_{\BQ}(X) := \NS(X) \otimes_{\BZ} \BQ$. Set
$$G^{*} = {\rm Im}(G\longrightarrow {\rm GL}\, (\NS_{\BQ}(X))).$$
Then the famous Tits Alternative Theorem says:

\begin{itemize}
\item [(i)]{\it 
either $G^{*}$ contains a subgroup isomorphic to
the non-abelian free group $\BZ * \BZ$ of rank two (highly noncommutative); or
}
\item[(ii)] {\it
$G^{*}$ contains a connected solvable subgroup $G_1$ of finite
index.
}
\end{itemize}

\noindent
Here, $G_1$ is {\it connected} if its Zariski closure in
${\rm GL}\, (\NS_{\BC}(X))$ is connected.

In \cite{DS04}, Dinh and Sibony proved the following
inspiring result:

\begin{theorem} [{\cite{DS04}}] \label{theorem:dsint}
Let $M$ be an $n$-dimensional compact K\"ahler manifold. Let
$G$ be an abelian subgroup of ${\rm Aut}\, (M)$ such that each
element of $G \setminus \{\id\}$ is of positive entropy. Then $G$
is a free abelian group of rank at most $n -1$. Furthermore, the
rank estimate is optimal.
\end{theorem}

In view of the Tits Alternative Theorem and Dinh-Sibony's Theorem, it is natural
to pose the following conjecture:

\begin{conjecture} \label{conjecture:titstype} {\bf (Conjecture of Tits type)}
\par \noindent
{\rm
Let $X$ be an $n$-dimensional compact K\"ahler manifold or an $n$-dimensional
complex projective variety with at most rational $\BQ$-factorial singularities.
Let $G$ be a subgroup of ${\rm Aut}\, (X)$. Then,
one of the following two assertions holds:
}
\begin{itemize}
\item[$(1)$] 
$G$ contains a subgroup isomorphic to $\BZ * \BZ$.
\item[$(2)$]
There is a finite-index subgroup $G_1$ of $G$ such that the
subset
$$N(G_1) = \{g \in G_1\, \vert\, g\,\ {\rm is\,\ of\,\ null\,\ entropy}\}\,$$
of $G_1$ is a normal subgroup of $G_1$ and the quotient group
$G_1/N(G_1)$ is a free abelian group of rank at most $n-1$.
\end{itemize}
\end{conjecture}

It turns out that the crucial part of the conjecture is the rank
estimate in the statement (2), as a weaker version of the
conjecture can be easily verified. 

\par
In \cite{KOZ}, we prove the result below. In the process, we need to
utilize Birkhoff's generalized Perron-Frobenius theorem and show
a Lie-Kolchin type theorem for nef cones; some geometrical analysis of
Calabi-Yau threefolds is also necessary.

\begin{theorem} [{cf. \cite{KOZ}}]
The conjecture of Tits type for a group $G$ on a compact complex variety $X$
holds in the following cases.
\begin{enumerate}
\item[$(1)$]
$\dim X \le 2$.
\item[$(2)$]
$X$ is a Hyperk\"ahler manifold.
\item[$(3)$]
$X$ is a complex torus.
\item[$(4)$]
$X$ is a projective minimal threefold
(with at worst $\BQ$-factorial terminal singularities
and nef canonical divisor).
\end{enumerate}
\end{theorem}

\section{Dynamics of automorphisms}
%
%
%
%

In this section, we show that the dynamics of automorphisms on all projective complex
manifolds (of dimension 3, or of any dimension but assuming the Good
Minimal Model Program or Mori's Program) are canonically built up
from the dynamics on just three types of projective complex
manifolds $X$ (here $X$ is $wCY$ if $\kappa(X) = 0 = q(X)$; see \cite{Zh2}, \cite{NZ}):
\par
{\it Complex Torus, wCY, Rationally
Connected Manifold.}

\begin{theorem}[{cf. \cite{Zh2}}]\label{Thk}
Let $X$ be a smooth projective complex manifold
of $\dim X \ge 2$, and with
$g \in \Aut(X)$. Then we have\emph:
\begin{itemize}
\item[$(1)$] Suppose that $(X, g)$ is either rigidly parabolic
or of primitively positive entropy. Then
the Kodaira dimension $\kappa(X) \le 0$.
\item[$(2)$]
Suppose that $\dim X = 3$ and $g$ is of positive entropy.
Then $\kappa(X) \le 0$, unless $d_1(g^{-1}) = d_1(g)= d_2(g) = e^{h(g)}$
and it is a Salem number. Here $d_i(g)$ are
dynamical degrees and $h(g)$ is the entropy.
\end{itemize}
\end{theorem}

\par
In the above, an algebraic integer $\lambda > 1$ of degree
$2(r+1)$ over $\BQ$ with $r \ge 0$, is a 
{\it Salem} number if all conjugates of $\lambda$ over $\BQ$ 
(including $\lambda$ itself) are given as follows, where
$|\alpha_i| = 1$:
$$\lambda, \, \lambda^{-1}, \, \alpha_1, \, {\bar \alpha}_1,\dots,\alpha_r,
\, \bar{\alpha}_r.$$

\par
In view of Theorem \ref{Thk}, we have only to treat the dynamics on
those $X$ with $\kappa(X) = 0$ or $-\infty$. 
This is done in \cite{Zh2}.
Below is an sample result for threefolds, which
says that $3$-dimensional dynamics of positive
entropy (not necessarily primitive) are just those of $3$-tori, wCY
$3$-folds and rational connected $3$-folds, unless
dynamical degrees are Salem numbers.
\begin{theorem}[{cf. \cite{Zh2}}]\label{3k}
Let $X'$ be a smooth projective complex threefold. Suppose that $g
\in \Aut(X')$ is of positive entropy.
Then there is a pair $(X, g)$ birationally equivariant to $(X', g)$,
such that one of the cases below occurs.
\begin{itemize}
\item[$(1)$]
There exist $3$-torus $\tilde{X}$ and
$g$-equivariant \'etale Galois cover $\tilde{X} \rightarrow X$.
\item[$(2)$]
$X$ is a wCY.
\item[$(3)$]
$X$ is a rationally connected threefold in the sense of Campana and Kollar-Miyaoka-Mori.
\item[$(4)$]
$d_1(g^{\pm} | X) = d_2(g^{\pm} | X) = e^{h(g | X)}$
and it is a Salem number.
\end{itemize}
\end{theorem}

We refer to McMullen \cite{Mc02a}, \cite{Mc02b}, \cite{Mc07} for the relation
between Salem numbers and K3 surfaces or anti-canonical rational surfaces.
See also McMullen \cite{Mc02b} and Cantat \cite{Ca} for the systematic study of 
dynamics on K3 surfaces.

\section{Building blocks of endomorphisms}
%
%
%
%

We will prove in this section the claim below.

\begin{claim} \label{claim}
{\rm
All \'etale
endomorphisms of projective manifolds $X$ are constructed from two building blocks below
(up to isomorphism), assuming the good minimal model conjecture (known up to dimension three):}

\par \vskip 1pc
Endomorphism of abelian varieties, and
\par
Nearly \'etale rational endomorphisms of weak Calabi-Yau varieties.

\end{claim}

\par
See Definitions \ref{wCY} and \ref{n-et} below for the notion above.
The \'etaleness assumption of $X$ is quite natural because every surjective 
endomorphism of $X$ is \'etale provided that $X$ is a nonsingular 
projective variety and is non-uniruled. 

\par
We like to use building blocks, via 'canonical fibration', to construct canonically all \'etale endomorphisms.
The canonicity here roughly corresponds to the equivariance of a fibration or morphism.
On a nonsingular projective variety $X$, there are three
such canonical fibrations:

\par \vskip 1pc
{\it The Iitaka fibration for those $X$ of Kodaira dimension $\kappa(X) > 0$,}

\par 
{\it The albanese map for those $X$ with irregularity $q(X) > 0$, and}

\par
{\it The maximal rationally connected fibration for uniruled varieties.}

\par \vskip 1pc
Our reduction of endomorphisms to the building blocks will go along the line of
these canonical fibrations. But we need to take care of the equivariance
which is not always true.

\par
We start with the result below for those $X$
of positive Kodaira dimension.

\par
Theorem \ref{ThA} below treats not only holomorphic surjective endomorphisms of projective 
varieties of $ \kappa > 0 $ but also 
dominant meromorphic endomorphisms of 
compact complex manifolds of $ \kappa > 0 $ 
in the class $ \SC $ in the sense of Fujiki.
Note that a compact complex manifold is in the class $ \SC $ if and 
only if it is bimeromorphic to a compact K\"ahler manifold.

\begin{theorem}[{cf. \cite[Theorem A]{NZ}}]\label{ThA}
Let $X$ be a compact complex manifold in the class $ \SC $ 
of $ \kappa(X) \geq 1 $ 
and let $f \colon X \ratmap X$ be a dominant meromorphic map. 
Let $ W_{m} $ be the image of the $ m $-th pluricanonical map 
\[\Phi_{m} \colon X \ratmap \LS{mK_{X}}^{\vee} 
= {\BP}(\OH^{0}(X, mK_{X}))\] 
giving rise to the Iitaka fibration of $ X $. 
Then there is an automorphism $g$ of $ W_{m} $ of 
finite order such that 
$\Phi_{m} \circ f = g \circ \Phi_{m}$. 
\end{theorem}

\begin{remark}\hfill 
{\rm
If $ f $ is holomorphic, then, resolving 
the indeterminacy points of $\Phi_{m}$, we may assume that both 
$f \colon X \to X$ and 
$\Phi_m \colon X \to W_m$ are holomorphic 
so that $\Phi_m \circ f = g\circ \Phi_m$. 
This is because $f$ is \'etale and we can take an equivariant 
resolution of the graph of Iitaka fibration.
%
}
\end{remark}

With Theorem \ref{ThD} below, we can see the reductive significance of 
Theorem \ref{ThA} above from the dynamics point of view.
Indeed, in Theorem \ref{ThA}, if we assume that $f$ is holomorphic then
$f$ is \'etale ($\kappa(X)$ being non-negative now);
we may assume that $g^s = \id$ (for some $s > 0$) so that $\Phi_m \circ f^s = \Phi_m$; 
thus the topological entropies and the first dynamical degrees satisfy 
$$(h_{\ttop}(f))^s = h_{\ttop}(f^s |_F), \hskip 2pc (d_1(f))^s = d_1(f^s|_F)$$
for a smooth fibre $F$ of $\Phi_m$
(see the remark above). 

\begin{theorem} [{cf. \cite[Appendix A: Th D, Th A.10]{NZ}}]\label{ThD}  $ $
\newline
Let $\pi : X \to Y$ be a fibre space from a compact K\"ahler manifold $X$ onto
a compact complex analytic variety $Y$ and let $f : X \to X$ be an \'etale surjective endomorphism
such that $\pi \circ f = \pi$. Then we have\emph:
\begin{itemize}
\item[$(1)$] The equality $h_{\ttop}(f) = h_{\ttop}(f|_F)$ holds
for the topological entropies $h_{\ttop}$ of $f$ and its restriction 
$f|_F : F \to F$ to a smooth fibre of $\pi$.
\item[$(2)$] The equality $d_1(f) = d_1(f|_F)$ holds
for the first dynamical degrees $d_1$ of $f$ and $f|_F$ (the \'etaleness of $f$
is not needed here).
\end{itemize}
\end{theorem}

\par
In view of Theorem \ref{ThA} above, we are reduced to considering endomorphisms $f : X \to X$
on those nonsingular projective varieties $X$ of Kodaira dimension zero 
(like the fibres $F$ of the Iitaka fibration in Theorem \ref{ThA}).

\par
Now let $X$ be a nonsingular projective variety of Kodaira dimension zero.
One would naturally try to reduce the study of $f : X \to X$ to that of the induced
self map on a minmal model $X_{\min}$ of $X$. 

\par
Indeed, the existence of a good minimal model
in the sense of Kawamata is known in dimension three or less.
In higher dimension, there has been rapid progress due to Birkar, Cascini, Hacon and McKernan \cite{BCHM}.

\par
The trouble is, on the one hand, that an holomorphic endomorphism $f : X \to X$ induces only 
a rational map $f_{X_{\min}} : X_{\min} \ratmap X_{\min}$.

\par
The good thing is, on the other hand, that $f_{X_{\min}}$ satisfies the condition of 
\emph{nearly \'etale map} to be defined below.
So our result below is applicable to $f_{X_{\min}}$.

\begin{theorem}[{cf. \cite[Theorem B]{NZ}}]\label{ThB}
Let $ V $ be a normal projective variety with only canonical 
singularities such that $ K_{V} \sim_{\BQ} 0 $. 
Let $ h \colon V$ $\ratmap V $ be a dominant rational map which is 
nearly \'etale. 
Then there exist an abelian variety $ A $, 
a weak Calabi--Yau variety $ F $ (see below), 
a finite \'etale morphism $ \tau \colon F \times A \to V $, 
a nearly \'etale dominant rational map 
$ \varphi_{F} \colon F \ratmap F $, and 
a finite \'etale morphism $ \varphi_{A} \colon A \to A $ 
such that 
$ \tau \circ (\varphi_{F} \times \varphi_{A}) = h \circ \tau $, 
i.e., the diagram below is commutative\emph{:}
\[ \begin{array}{ccc}
F \times A & \overset{\varphi_{F} \times \varphi_{A}}{\ratmap} & 
F \times A \\
\mbox{\scriptsize $\tau$} \downarrow \phantom{\tau} & & \phantom{\tau} 
\downarrow \mbox{\scriptsize $\tau$} \\
V &\overset{h}{\ratmap} & \phantom{.}V.
\end{array}  \]
\end{theorem}

\begin{remark}\hfill
{\rm 
\begin{enumerate}

\item[$(1)$]  
If $ F $ has only terminal singularities, then $ \varphi_{F} $ 
is \'etale in codimension one. 

\item[$(2)$] 
If the algebraic fundamental group $ \pi_{1}^{\alg}(F)$ is finite
(this is true if $\dim F \le 3$ by \cite[Corollary (1.4)]{NS}),
then $ \varphi_{F} $ is a birational automorphism.
In particular, if $ V $ has only terminal singularities and $ 
q^{\max}(V) = \dim V - 2 $,
then $ \varphi_{F} $ is an automorphism.  
\end{enumerate}
}
\end{remark}

\begin{definition}\label{wCY}
{\rm
A normal projective variety $F$ is called 
\emph{weak Calabi--Yau} if 
$ F $ has only canonical singularities, $ K_{F} \sim_{\BQ} 0 $, and 
\[ q^{\max}(F) := \max\{ q(F') \mid F' \to F \text{ is finite \'etale} 
\} = 0 . \] 
\par
If $ F $ is a nonsingular weak Calabi--Yau variety, 
then $ \pi_1(F) $ is finite by 
Bogomolov's decomposition theorem, so a finite \'etale cover of $F $ is expressed as 
a product of holomorphic symplectic manifolds 
and of Calabi--Yau manifolds. If $F$ is singular,
we have to use a result of Kawamata instead.
}
\end{definition}

\begin{definition}\label{n-et} 
{\rm
Let $ h \colon V \ratmap W $ be a proper 
rational (resp.\ meromorphic) map between 
algebraic (resp.\ complex analytic) varieties. 
The map $ h $ is called \emph{nearly \'etale} if 
there exist proper birational (resp.\ bimeromorphic) maps 
$ \mu \colon Y \ratmap W $, $ \nu \colon X \ratmap V $ 
and a morphism $ f \colon X \to Y $ such that 
\begin{enumerate}
\item[$(1)$]  
$X$ and $Y$ are algebraic (resp.\ complex analytic) varieties, 

\item[$(2)$]  
$ f $ is a finite \'etale morphism, and

\item[$(3)$]  
$ \mu \circ f = h \circ \nu $. 
\end{enumerate}
}
\end{definition}

As we have seen from the results above, the building blocks of
surjective endomorphisms of nonsingular projective varieties
with Kodaira dimension $\ge 0$, are the endomorphisms of abelian varieties
and nearly \'etale rational endomorphisms of weak Calabi-Yau varieties.

\par
We still have to treat the case of Kodaira dimension $-\infty$.
Let $X$ be a nonsingular projective variety of Kodaira dimension $-\infty$
and $f : X \to X$ a surjective endomorphism. 
The study of such $f$ is very difficult. For instance, we still
do not know well the structure of such $f$ even when $X$ is a {\it rational} surface,
though for surfaces $X$, the geometrical structure of the underlying
surface $X$ is well understood (provided that $f$ is non-isomorphic),
thanks to \cite{Na} and \cite{FN05}.

\par
On the other hand, it is conjectured that being $\kappa(X) = -\infty$
is equivalent to the uniruledness of $X$ ({\it weak Abundance conjecture}).
This conjecture is known if $\dim X \le 3$.

\par
Therefore, the result below (which uses \cite{BCHM})
and the weak Abundance conjecture reduce the study of 
$f : X \to X$ to that of an \'etale endomorphism $f_Y : Y \to Y$
of a non-uniruled projective manifold $Y$.
\par
Indeed, the $M$ in Theorem \ref{ThC} below is birational to
the Cartesian product $M' := M \times_Y Y$ 
of $\pi : M \to Y$ and $f_Y : Y \to Y$,
and $(X, f)$ and $(M, f_M)$ are birationally equivariant 
to $(M', \ {\rm pr}_M)$ with ${\rm pr}_M : M' \to M$ the first projection
(so $f$ and $f_M$ are birationally determined by $f_Y$);
this is because rationally connected projective manifolds
are simply connected.

\par
For $h : Y \to Y$ we could apply Theorem \ref{ThA}
(assuming weak Abundance conjecture), and note that $\dim Y < \dim X$.
This way, we have confirmed Claim \ref{claim}.

\begin{theorem}[{cf. \cite[Theorem C and its Remark]{NZ}}]\label{ThC}
$ $
\newline
Let $X$ be a projective manifold with an \'etale endomorphism $f$. 
Assume that $ X $ is uniruled. 
Lifting $(X, f)$ birationally equivariantly to a pair
$(M, f_M)$ (with $f_M : M \to M$ \'etale too), we have a maximal rationally connected fibration
$\pi : M \to Y$ onto a non-uniruled projective 
manifold $Y$ and an \'etale endomorphism $f_Y : Y \to Y$
such that $\pi \circ f_M = f_Y \circ \pi$,
i.e., the diagram below is commutative\emph{:}
\[ \begin{array}{ccc}
M & \overset{f_M}{\longrightarrow} & M \\
\mbox{\scriptsize $\pi$} \downarrow \phantom{\pi} & & \phantom{\pi} 
\downarrow \mbox{\scriptsize $\pi$} \\
Y &\overset{f_Y}{\longrightarrow} & \phantom{.}Y.
\end{array}  \]
\end{theorem}

\par
We refer to \cite{Fu} and \cite{FN07} for the precise geometrical structure of 
nonsingular projective threefolds with non-negative Kodaira dimension and admitting 
non-isomorphic surjective endomorphisms.
\par
We refer to \cite{Zh2} or \S 5 for the building blocks of automorphisms 
of algebraic manifolds.

\section{Cohomologically hyperbolic endomorphism}

%
%
%
%

\par
Let $X$ be a compact K\"ahler manifold of dimension $n \ge 2$.
A surjective endomorphism $f : X \to X$ is {\it cohomologically hyperbolic} 
in the sense of \cite{Gu06},
if there is an $\ell \in \{1, 2,  \dots, n\}$
such that the $\ell$-th dynamical degree 
$$d_{\ell}(f) > d_i(f) \hskip 2pc
{\rm for \hskip 1pc all}  \hskip 1pc (\ell \ne) \,\, i \,\, \in \,\, \{0, 1, \dots, n\},$$ 
or equivalently, for both $i = \ell \pm 1$, by
the Khovanskii - Tessier inequality. 
%
%
\par
In his papers \cite{Gu05} - \cite{Gu06}, Guedj assumed that a dominant rational self map $f : X \ratmap X$ 
has {\it large topological degree} (i.e., it is cohomologically hyperbolic with $\ell = \dim X$
in the definition above), and constructed a {\it unique} $f_*$-invariant measure $\mu_f$. Further,
the measure is proved to be of maximal entropy, ergodic,
equidistributive for $f$-periodic and repulsive points, and with strictly
positive Lyapunov exponents.
In \cite{Gu06}, Guedj classified cohomologically
hyperbolic rational self maps of surfaces $S$ and deduced that the Kodaira dimension
$\kappa(S) \le 0$. {\it Then he conjectured that the same should hold in higher dimension.}
\par
The result below gives an affirmative answer
to the above-mentioned conjecture of Guedj \cite{Gu06} page 7
for holomorphic endomorphisms (see \cite[Theorem 1.3]{Zh2} for the case
of automorphisms on threefolds).
The proof is given by making use of results in \cite{NZ} or Theorem \ref{ThD} above.
It is classification-free and for arbitrary dimension.

\begin{theorem} [{cf. \cite{Zh3}}]
Let $X$ be a compact complex K\"ahler manifold and
$f : X \to X$ a surjective and cohomologically hyperbolic endomorphism. 
Then the Kodaira dimension
$\kappa(X) \le 0$.
\end{theorem}

\newpage

\par
We now determine the geometric structure for projective threefolds in
the above theorem.

\begin{theorem}[{cf. \cite{Zh3}}]
Let $V$ be a smooth projective threefold and let $f \colon V \to V$ be a
surjective and cohomologically hyperbolic \'etale endomorphism.
Then one of the following cases occurs.
\begin{itemize}
\item[$(1)$]
$V$ is $f$-equivariantly birational to a {\rm Q}-torus,
i.e., the quotient of an abelian variety modulo a finite and free action.
\item[$(2)$]
$V$ is birational to a weak Calabi-Yau variety, and $f \in \Aut(V)$.
\item[$(3)$]
$V$ is rationally connected, 
and $f \in \Aut(V)$.
\item[$(4)$]
The albanese map $V \to \Alb(V)$ is a smooth and surjective morphism onto 
the elliptic curve $\Alb(V)$ with every fibre
a smooth projective rational surface of Picard number $\ge 11$.
Further, the dynamical degrees
satisfy $d_2(f) > d_1(f) \ge \deg(f) \ge 2$.
\item[$(5)$]
$V$ is $f$-equivariantly birational to the quotient space
of a product ${\rm (Elliptic \ curve)} \times$ $(K3)$ modulo a finite 
and free action. Further, the dynamical degrees satisfy
$d_2(f) > d_1(f) \ge \deg(f) \ge 2$.
\end{itemize}
\end{theorem}

\par \vskip 1pc \noindent
{\bf ACKNOWLEDGMENT}

\par \vskip 0.5pc
I would like to thank Noboru Nakayama for the constructive 
discussion and giving plenty of suggestions,
and the following colleagues for the discussion and comments:
\par \noindent
F. Campana, S. Cantat, F. Catanese, M. Chen,
T. -C. Dinh, H. Esnault, A. Fujiki, Y. Fujimoto, J. -M. Hwang,
Y. Kawamata,  C. Keem, 
J. Keum, S. Kondo, F. Kutzschebauch, N. C. Leung, W. -P. Li, S. Mukai,
K. Oguiso, N. Sibony, S. -L. Tan, E. Viehweg, S. -W. Zhang and K. Zuo.
\par
I also like to thank the following institutes for the support and warm
hospitality during the years 2006-2007:
\par \noindent
National Seoul Univ., Korea Inst. Adv. Study, Hong Kong Univ. Sci. Tech.
Chinese Univ. Hong Kong, 
Math. Forschungsinst. Oberwolfach,
Mainz Univ., Bayreuth Univ.,
Max-Planck-Inst. Math. at Bonn,
Essen Univ., Bern Univ., 
Univ. Paris 11-Orsay, Nagoya Univ.,
Univ. Tokyo, Osaka Univ. and
Res. Inst. Math. Sci. Kyoto Univ.
\par
The author wrote the preprint \cite{Zh2} on the building blocks of automorphisms in Autumn 2006;
a timely visit of Noboru Nakayama in October 2006 resulted in the paper \cite{NZ}
on the building blocks of endomorphisms. We would like to 
thank the Academic Research Fund of National University of Singapore, which made the visit possible.

\end{large}
\end{document}